\documentclass{article}

\newtheorem{thm}{Theorem}
\newtheorem{cor}[thm]{Corollary}
\newtheorem{lem}[thm]{Lemma}

\begin{document}

\title{ The probability that random positive integers are $k$-wise relatively prime}
\author{ Jerry Hu}
\date{}
\maketitle

\begin{abstract}
\hspace{-0.1957in}The positive integers $a_1,a_2, ...,a_s$ are $k$-wise relatively prime if any $k$ of them are relatively prime. Exact formula is obtained for the probability that s positive integers are $k$-wise relatively prime.
\end{abstract}

A classic result in number theory is that the probability that two given integers are relatively prime is $6/\pi^2$.
More generally the probability that $s$ positive integers chosen arbitrarily and independently are relatively prime
is $ 1/\zeta(s)$, where $\zeta(s)$ is Riemann's zeta function. A short accessible proof of this result is in the paper of Nymann \cite{nymann}. Recently, in 2002 T\'{o}th \cite{toth} solved the problem of finding the probability that $s$ positive integers are pairwise relatively prime by the recursion method.

The positive integers $a_1,a_2, ...,a_s$ are $k$-wise relatively prime if any $k$ of them are relatively prime and are $k$-wise relatively prime to $a$ if gcd of any $k$ of them is relatively prime to $a$. In this note we consider the problem of finding the probability $A_{s,k}$ that s positive integers are $k$-wise relatively prime.

For a $(k-1)$-tuple of positive integers $u = (u_1, ...,u_{k-1})$,
let $ Q_{s,k}^{(u)}(n)$ denote the number of $s$-tuples of positive integers $a_1,a_2, ...,a_s $ with $1 \leq a_1,a_2,..., a_s \leq n$ such that $a_1, a_2,...,a_s$ are $k$-wise relatively prime and are $i$-wise relatively prime to $u_i$ for $i = 1, 2, ..., k-1$.

The next theorem gives an asymptotic formula for $Q_{s,k}^{(u)}(n)$ and the exact values of $A_{s,k}$.
\begin{thm} For fixed $s\geq 1$, $k\geq 2$,  we have uniformly for $n, u_i\geq 1$ with $(u_i,u_j)=1$ for $1 \leq i\not = j \leq k-1,$
\begin{equation}  Q_{s,k}^{(u)}(n)=A_{s,k}\prod_{i=1}^{k-1}f_{s,k,i}(u_i)n^s +  O(\theta(u_1)n^{s-1}\log^{\delta(s,k)}\hspace{-0.05in} n ),  \label{meqn}
\end{equation}
where

$$A_{s,k}=\prod_p(1-\frac{1}{p})^{s-k+1}\sum_{m=0}^{k-1} \left(\hspace{-0.07in}\begin{array}{c} s \\ m  \end{array}\hspace{-0.07in}\right) (1-\frac{1}{p})^{k-1-m}\frac{1}{p^m},$$
$$f_{s,k,i}(u_i)=\prod_{p|u_i}\left(1-\frac{\sum_{m=i}^{k-1}\left(\hspace{-0.07in}\begin{array}{c} s \\ m  \end{array}\hspace{-0.07in}\right)(p-1)^{k-1-m}}{\sum_{m=0}^{k-1}\left(\hspace{-0.07in}\begin{array}{c} s \\ m  \end{array}\hspace{-0.07in}\right)(p-1)^{k-1-m}}\right),$$
$$ \delta(s,k)=\max\left\{\left(\hspace{-0.08in}\begin{array}{c} s-1 \\ i  \end{array}\hspace{-0.08in}\right), i=1,2,...,k-1 \right\}, $$

and $\theta (u_1)$ is the number of squarefree divisors of $u_1$.
\end{thm}

Note that $f_{s,k,i}(u_i)$ is a multiplicative arithmetic function for $i=1,...,k-1.$
\begin{cor}
The probability that s positive integers $a_1, a_2, ..., a_s$ are $k$-wise relatively prime and are $i$-wise relatively prime to $u_i$ for $i=1,2,...,k-1$ is
$$ \lim_{n \rightarrow \infty} \frac{Q_{s,k}^{(u)}(n)}{n^s}=A_{s,k} \prod_{i=1}^{k-1}f_{s,k,i}(u_i).$$
For $u_1=u_2=\cdot\cdot\cdot=u_{k-1}=1$, the probability that s positive integers are $k$-wise relatively prime is
$$  A_{s,k}=\prod_p(1-\frac{1}{p})^{s-k+1}\sum_{m=0}^{k-1} \left(\hspace{-0.07in}\begin{array}{c} s \\ m  \end{array}\hspace{-0.07in}\right) (1-\frac{1}{p})^{k-1-m}\frac{1}{p^m}.$$
\end{cor}

For a pair of positive integers $a$ and $b$, let $(a, b]$ denote the product of prime powers $p^e$ such that $p^e$ divides $b$, but $p^{e+1}$ does not divide $b$, and $p$ also divides $a$. Note that $(a, b)$ and $b/(a,b)$ may not be relatively prime, but $(a,b)$ and $b/(a,b]$ are relatively prime.

To prove the theorem we need the following lemmas.

\begin{lem}
For $ s,n \geq 1$, $k\geq 2$, and $ u=(u_1,u_2,...,u_{k-1})$ with $ u_1,...,u_{k-1} \geq 1$ and $ (u_i, u_j)=1$ for $ i \neq j$,
$$ Q_{s+1,k}^{(u)}(n)  =  \sum_{\scriptsize \begin{array}{c} j=1\\(j,u_1)=1 \end{array}}^nQ_{s,k}^{(j\ast u)}(n),$$\
where
$$ j* u=\left(u_1(j,u_2), \frac{u_2(j,u_3)}{(j,u_2]}, ...,\frac{u_{k-2}(j,u_{k-1})}{(j,u_{k-2}]},\frac{ju_{k-1}}{(\prod_{i=2}^{k-1}[j,u_i))(j,u_{k-1}]}\right).$$
\end{lem}
{\em Proof.} An $(s+1)$-tuple of positive integers $a_1,a_2,...,a_{s+1}$ are $k$-wise relatively prime and are $i$-wise relatively prime to $u_i$ for $i=1,2,...,k-1$  if and only if the first $s$-tuple of positive integers $a_1,a_2,...,a_s$ are $k$-wise relatively prime and are $i$-wise relatively prime to $u_i$ and $(a_{s+1},u_{i+1})$ for $i=1,2,...,k-2$ and are $(k-1)$-wise relatively prime to $u_{k-1}$ and $a_{s+1}$, and $(a_{s+1}, u_1)$=1,  we have
$$ Q_{s+1,k}^{(u)}(n)= \sum_{\scriptsize \begin{array}{c} a_{s+1}=1\\(a_{s+1},u_1)=1 \end{array}}^nQ_{s,k}^{(a_{s+1}*' u)}(n)=\sum_{\scriptsize \begin{array}{c} j=1\\(j,u_1)=1 \end{array}}^nQ_{s,k}^{(j*' u)}(n),$$
where
$$ j*' u = (u_1(j,u_2), u_2(j,u_3),...,u_{k-2}(j,u_{k-1}),ju_{k-1}).$$

Note that the $(k-1)$-tuples of positive integers $j*' u$ are not pairwise relatively prime, so Theorem 1 can not be applied to $Q_{s,k}^{(j*' u)}(n)$ and the above recurrence relation can not be used in the proof of Theorem 1 by induction.

 To complete the proof of the lemma, it suffices to show that $Q_{s,k}^{(j*' u)}(n)=Q_{s,k}^{(j* u)}(n)$. For an $s$-tuple of positive integers, if each of them is relatively prime to $(j,u_2)$, then they are pairwise relatively prime to $(j,u_2),$ and pairwise relatively prime to $(j,u_2]$ since $(j,u_2)$ and $(j,u_2]$ contain the same prime factors, thus they are pairwise relatively prime to $u_2(j,u_3)$ if and only if they are pairwise relatively prime to $u_2(j,u_3)/(j,u_2]$. Similarly, for $i=3,4,...,k-2,$ if an $s$-tuple of positive integers are $(i-1)$-wise relatively prime to $(j,u_i)$, then they are $i$-wise relatively prime to $u_i(j,u_{i+1})$ if and only if they are $i$-wise relatively prime to $u_i(j,u_{i+1})/(j,u_i].$ For $i=1,2,...,k-2,$ if an $s$-tuple of positive integers are $i$-wise relatively prime to $(j,u_{i+1})$, then they are $(k-1)$-wise relatively prime to $[j,u_{i+1})$, and $(k-1)$-wise relatively prime to $(j,u_{k-1}]$, thus they are $(k-1)$-wise relatively prime to $ju_{k-1}$ if and only if they are $(k-1)$-wise relatively prime to $ju_{k-1}/(\prod_{i=2}^{k-1}[j,u_i))(j,u_{k-1}].$ Therefore, $Q_{s,k}^{(j*' u)}(n)=Q_{s,k}^{(j* u)}(n). $

Note that the $(k-1)$-tuple of integers in $j*u$ are pairwise relatively prime, so we can apply Theorem 1 to $Q_{s,k}^{(j* u)}(n)$ in the proof of Theorem 1 by induction.

\begin{lem}
For $k, u_i \geq 1$,
$$ \frac{f_{s,k,i}(u_i)}{f_{s,k,i+1}(u_i)} = \sum_{d | u_i}\frac{\mu(d)\left(\hspace{-0.07in}\begin{array}{c} s \\ i \end{array}\hspace{-0.07in}\right)^{\omega(d)}}{\alpha_{s,k,i}(d)}, i=1,2,...,k-2,$$
$$ f_{s,k,k-1}(u_{k-1})= \sum_{d|u_{k-1}}\frac{\mu(d)\left(\hspace{-0.08in}\begin{array}{c} s \\ k-1  \end{array}\hspace{-0.08in}\right)^{\omega(d)}}{\alpha_{s,k,k-1}(d)},$$
where
$$ \alpha_{s,k,i}(d) = d^i \prod_{p |d}\sum_{m=0}^{i} \left(\hspace{-0.07in}\begin{array}{c} s \\ m \end{array}\hspace{-0.07in}\right) (1-\frac{1}{p})^{i-m}\frac{1}{p^m},i=1,2,...,k-1,$$
and $\omega(u_i)$ denote the number of distinct prime factors of $u_i.$
\end{lem}
{\em Proof.} Since both $\frac{f_{s,k,i}(u_i)}{f_{s,k,i+1}(u_i)}$ and $f_{s,k,k-1}(u_{k-1})$ are multiplicative arithmetic functions, it suffices to  verify for $u_i=p^a$ a prime power:

\begin{eqnarray*}  \sum_{d | p^a}\frac{\mu(d)\left(\hspace{-0.07in}\begin{array}{c} s \\ i \end{array}\hspace{-0.07in}\right)^{\omega(d)}}{\alpha_{s,k,i}(d)}
& = & 1-\left( \hspace{-0.07in}\begin{array}{c} s \\ i \end{array}\hspace{-0.07in}\right) \left(\sum_{m=0}^{i} \left(\hspace{-0.07in}\begin{array}{c} s \\ m  \end{array}\hspace{-0.07in}\right) (p-1)^{i-m}\right)^{\hspace{-0.07in}-1} = \frac{f_{s,k,i}(p^a)}{f_{s,k,i+1}(p^a)},\\
\sum_{d|p^a}\frac{\mu(d)\left(\hspace{-0.08in}\begin{array}{c} s \\ k-1  \end{array}\hspace{-0.08in}\right)^{\omega(d)}}{\alpha_{s,k,k-1}(d)}& = &  1-\left(\hspace{-0.08in}\begin{array}{c} s \\ k - 1  \end{array}\hspace{-0.08in}\right)\left( \sum_{m=0}^{k-1} \left(\hspace{-0.07in}\begin{array}{c} s \\ m  \end{array}\hspace{-0.07in}\right) (p-1)^{k-1-m} \right)^{\hspace{-0.07in}-1} = f_{s,k,k-1}(p^a).
\end{eqnarray*}

For the proof of the theorem, we proceed by induction on $s$. For $s=1, $ we have by the Inclusion-Exclusion Principle
\begin{eqnarray*}
Q_{1, k}^{(u)}(n) & = & \sum_{\scriptsize \begin{array}{c} j=1 \\(j,u_1)=1 \end{array}}^{n} 1=\sum_{d|u_1}\mu(d)\lfloor\frac{n}{d}\rfloor=\sum_{d|u_1}\mu(d)\left(\frac{n}{d}+O(1)\right)\\
 & = & n\sum_{d|u_1}\frac{\mu(d)}{d}+O\left(\sum_{d|u_1}\mu^2(d)\right).
\end{eqnarray*}
Hence,
\begin{equation} Q_{1,k}^{(u)}(n)  =  \sum_{\scriptsize \begin{array}{c} j=1 \\(j,u_1)=1 \end{array}}^{n}1=n\frac{\phi(u_1)}{u_1}+O(\theta(u_1)) \label{eueqn} \end{equation}
and (\ref{meqn}) is true for $s=1$ with $A_{1,k}=1, f_{1,k,i}(u_i)=1$ for $i \geq 2$, and $f_{1,k,1}(u_1)=\frac{\phi(u_1)}{u_1}, \phi$ denoting the Euler function.

Suppose that (\ref{meqn}) is valid for $s$, we prove it for $s+1$. From Lemma $3$, we have

\begin{eqnarray}
Q_{s+1,k}^{(u)}(n) & = & \sum_{\scriptsize \begin{array}{c} j=1\\(j,u_1)=1 \end{array}}^nQ_{s,k}^{(j*u)}(n)\nonumber \\
& = & \sum_{\scriptsize \begin{array}{c} j=1\\(j,u_1)=1 \end{array}}^n A_{s,k} \prod_{i=1}^{k-2}f_{s,k,i}\left(\frac{u_i(j,u_{i+1})}{(j,u_i]}\right)f_{s,k,k-1}\left(\frac{ju_{k-1}}{(\prod_{i=2}^{k-1}[j,u_i))(j,u_{k-1}]}\right)n^s \nonumber \\ & & \hspace{1.1in}  \mbox{} +O(\theta(u_1(j,u_2))n^{s-1}\log^{\delta(s,k)}\hspace{-0.05in}n) \label{eqnb}  \\
 & = & A_{s,k}\prod_{i=1}^{k-1} f_{s,k,i}(u_i)n^s \sum_{\scriptsize \begin{array}{c} j=1\\(j,u_1)=1 \end{array}}^n \prod_{i=1}^{k-2}\frac{f_{s,k,i}((j,u_{i+1}))}{f_{s,k,i+1}((j,u_{i+1}))}f_{s,k,k-1}\left(\frac{j}{\prod_{i=2}^{k-1}[j,u_i)}\right) \nonumber \\
 & & \hspace{1.1in} \mbox{}+ O(\theta(u_1)n^{s-1}\log^{\delta(s,k)}\hspace{-0.05in}n\sum_{j=1}^n\theta(j)). \nonumber
\end{eqnarray}
Here $\sum_{j=1}^n\theta(j) \leq \sum_{j=1}^n\tau_2(j)=O(n\log n),$ where $\tau_2 = \tau$ is the divisor \mbox{function}.

Furthermore, from Lemma 4,
\begin{eqnarray*}
\lefteqn{\sum_{\scriptsize \begin{array}{c} j=1\\(j,u_1)=1 \end{array}}^n \prod_{i=1}^{k-2}\frac{f_{s,k,i}((j,u_{i+1}))}{f_{s,k,i+1}((j,u_{i+1}))}f_{s,k,k-1}\left(\frac{j}{\prod_{i=2}^{k-1}[j,u_i)}\right)}\\
& = & \sum_{\scriptsize \begin{array}{c} j=1\\(j,u_1)=1 \end{array}}^n\prod_{i=1}^{k-2}\sum_{d_i\mid (j,u_{i+1})}\frac{\mu(d_i)\left(\hspace{-0.07in}\begin{array}{c} s \\ i  \end{array}\hspace{-0.07in}\right)^{\omega(d_i)}}{\alpha_{s,k,i}(d_i)}\sum_{d_{k-1} |\frac{j}{\prod_{i=2}^{k-1}[j,u_i)}}\frac{\mu(d_{k-1})\left(\hspace{-0.08in}\begin{array}{c} s \\ k-1  \end{array}\hspace{-0.08in}\right)^{\omega(d_{k-1})}}{\alpha_{s,k,k-1}(d_{k-1})}\\
 & = & \sum_{\scriptsize \begin{array}{c} d_1d_2\cdot\cdot\cdot d_{k-1}e=j\leq n\\d_i | (j,u_{i+1}), i=1,...,k-2\\d_{k-1} | \frac{j}{[j,u_2)\cdot\cdot\cdot[j,u_{k-1})}\\ (j,u_1)=1 \end{array}} \prod_{i=1}^{k-1}\frac{\mu(d_i)\left(\hspace{-0.07in}\begin{array}{c} s \\ i  \end{array}\hspace{-0.07in}\right)^{\omega(d_i)}}{\alpha_{s,k,i}(d_i)}\\
 &  = & \sum_{\scriptsize \begin{array}{c} d_1d_2\cdot\cdot\cdot d_{k-1}\leq n\\d_i | u_{i+1}, i=1,...,k-2\\(d_{k-1},u_i)=1, i=1,2,...,k-1 \\ \end{array}} \prod_{i=1}^{k-1}\frac{\mu(d_i)\left(\hspace{-0.07in}\begin{array}{c} s \\ i  \end{array}\hspace{-0.07in}\right)^{\omega(d_i)}}{\alpha_{s,k,i}(d_i)}
\sum_{\scriptsize \begin{array}{c} e\leq\frac{n}{d_1d_2\cdot\cdot\cdot d_{k-1}}\\(e,u_1)=1\end{array}}1
\end{eqnarray*}

Using  (\ref{eueqn}), we have
\begin{eqnarray}
\lefteqn{\sum_{\scriptsize \begin{array}{c} j=1\\(j,u_1)=1 \end{array}}^n \prod_{i=1}^{k-2}\frac{f_{s,k,i}((j,u_{i+1}))}{f_{s,k,i+1}((j,u_{i+1}))}f_{s,k,k-1}\left(\frac{j}{\prod_{i=2}^{k-1}[j,u_i)}\right)}\\
  & =& \sum_{\scriptsize \begin{array}{c} d_1d_2\cdot\cdot\cdot d_{k-1}=d \leq n\\d_i | u_{i+1}, i=1,...,k-2\\(d_{k-1},u_i)=1, i=1,2,...,k-1 \\ \end{array}} \prod_{i=1}^{k-1}\frac{\mu(d_i)\left(\hspace{-0.07in}\begin{array}{c} s \\ i  \end{array}\hspace{-0.07in}\right)^{\omega(d_i)}}{\alpha_{s,k,i}(d_i)}
 \left(\frac{\phi(u_1)}{u_1}\frac{n}{d_1d_2\cdot\cdot\cdot d_{k-1}}+O(\theta(u_1))\right) \nonumber \\
 & =& \frac{\phi(u_1)}{u_1} n \sum_{\scriptsize \begin{array}{c} d_1d_2\cdot\cdot\cdot d_{k-1}=d \leq n\\d_i | u_{i+1}, i=1,...,k-2\\(d_{k-1},u_i)=1, i=1,2,...,k-1 \\ \end{array}} \prod_{i=1}^{k-1}\frac{\mu(d_i)\left(\hspace{-0.07in}\begin{array}{c} s \\ i  \end{array}\hspace{-0.07in}\right)^{\omega(d_i)}}{d_i\alpha_{s,k,i}(d_i)}
  \label{eqna} \\
 & & \hspace{1.1in} \mbox{}+ O\left(\theta(u_1)\sum_{\scriptsize d \leq n}                            \frac{\delta(s+1,k)^{\omega(d)}}{d}\right), \nonumber
\end{eqnarray}
since $\alpha_{s,k,i}(d_i) > d_i.$

Hence, the main term  of (\ref{eqna}) is
\begin{eqnarray*}
\lefteqn{\frac{\phi(u_1)}{u_1} n \sum_{\scriptsize \begin{array}{c} d_i | u_{i+1}, i=1,...,k-2\\(d_{k-1},u_i)=1, i=1,2,...,k-1 \\ \end{array}} \prod_{i=1}^{k-1}\frac{\mu(d_i)\left(\hspace{-0.07in}\begin{array}{c} s \\ i  \end{array}\hspace{-0.07in}\right)^{\omega(d_i)}}{d_i\alpha_{s,k,i}(d_i)}}\\
 &=& \hspace{-0.1in} \frac{\phi(u_1)}{u_1} n \hspace{-0.03in}\prod_{i=2}^{k-1}\prod_{p | u_i}\hspace{-0.04in}\left(1-\frac{\left(\hspace{-0.08in}\begin{array}{c} s \\ i-1 \end{array}\hspace{-0.08in}\right)}{p\sum_{m=0}^{i-1}\left(\hspace{-0.07in}\begin{array}{c} s \\ m  \end{array}\hspace{-0.07in}\right)(p-1)^{i-1-m}}\right)\hspace{-0.07in} \prod_{ p\hspace{-0.01in} \not \hspace{0.04in}|  u_{\hspace{-0.02in}1}\cdot\cdot\cdot u_{k\hspace{-0.015in}-\hspace{-0.025in}1}}\hspace{-0.07in} \left(1-\frac{\left(\hspace{-0.08in}\begin{array}{c} s \\ k-1  \end{array}\hspace{-0.08in}\right)}{p\sum_{m=0}^{k-1}\left(\hspace{-0.07in}\begin{array}{c} s \\ m  \end{array}\hspace{-0.07in}\right)(p-1)^{k-1-m}}\right) \\
 &=& \hspace{-0.1in}n \prod_{i=1}^{k-1}\prod_{p | u_i}\left(1-\frac{\left(\hspace{-0.08in}\begin{array}{c} s \\ i-1 \end{array}\hspace{-0.08in}\right)}{p\sum_{m=0}^{i-1}\left(\hspace{-0.07in}\begin{array}{c} s \\ m  \end{array}\hspace{-0.07in}\right)(p-1)^{i-1-m}}\right)
 \left(1-\frac{\left(\hspace{-0.08in}\begin{array}{c} s \\ k-1  \end{array}\hspace{-0.08in}\right)}{p\sum_{m=0}^{k-1}\left(\hspace{-0.07in}\begin{array}{c} s \\ m  \end{array}\hspace{-0.07in}\right)(p-1)^{k-1-m}}\right)^{-1} \\
 & & \prod_{p }\left(1-\frac{\left(\hspace{-0.08in}\begin{array}{c} s \\ k-1  \end{array}\hspace{-0.08in}\right)}{p\sum_{m=0}^{k-1}\left(\hspace{-0.07in}\begin{array}{c} s \\ m  \end{array}\hspace{-0.07in}\right)(p-1)^{k-1-m}}\right),
 \end{eqnarray*}

 and its O-terms are
 \begin{eqnarray*}
 O\left(n\sum_{d > n}\frac{\delta(s+1,k)^{\omega(d)}}{d^2}\right ) &=&       O\left(n\sum_{d > n}\frac{\tau_{\delta(s+1,k)}(d)}{d^2}\right) \\
     & = & O(\log^{\delta(s+1,k)-1}\hspace{-0.05in}n)
 \end{eqnarray*}
by Lemma 3(b) in \cite{toth}, which gives an asymptotic estimate of the sum $$ \sum_{n > x}\frac{\tau_k(n)}{n^2}=O\left(\frac{\log^{k-1}x}{x}\right)$$
  and

\begin{eqnarray*}
 O\left(\theta(u_1)\sum_{d \leq n}\frac{\delta(s+1,k)^{\omega(d)}}{d}\right)
& = & O\left(\theta(u_1)\sum_{d \leq n}\frac{\tau_{\delta(s+1,k)}(d)}{d}\right)\\
&=& O(\theta(u_1)\log^{\delta(s+1,k)}\hspace{-0.05in}n)
\end{eqnarray*}
from Lemma 3(a) in \cite{toth}, which gives an asymptotic estimate of the sum $$ \sum_{n\leq x}\frac{\tau_k(n)}{n}=O(\log^{k}x). $$

Substituting into (\ref{eqnb}), we get
\begin{eqnarray*}
Q_{s+1,k}^{(u)}(n) & = &A_{s,k}\prod_p\left(1-\frac{\left(\hspace{-0.08in}\begin{array}{c} s \\ k-1  \end{array}\hspace{-0.08in}\right)}{p\sum_{m=0}^{k-1}\left(\hspace{-0.07in}\begin{array}{c} s \\ m  \end{array}\hspace{-0.07in}\right)(p-1)^{k-1-m}} \right)\\
 \lefteqn{\hspace{-0.2in} \prod_{i=1}^{k-1}f_{s,k,i}\prod_{p | u_i}\left(1-\frac{\left(\hspace{-0.08in}\begin{array}{c} s \\ i-1 \end{array}\hspace{-0.08in}\right)}{p\sum_{m=0}^{i-1}\left(\hspace{-0.07in}\begin{array}{c} s \\ m  \end{array}\hspace{-0.07in}\right)(p-1)^{i-1-m}}\right)
 \left(1-\frac{\left(\hspace{-0.08in}\begin{array}{c} s \\ k-1  \end{array}\hspace{-0.08in}\right)}{p\sum_{m=0}^{k-1}\left(\hspace{-0.07in}\begin{array}{c} s \\ m  \end{array}\hspace{-0.07in}\right)(p-1)^{k-1-m}}\right)^{\hspace{-0.08in}-1}\hspace{-0.08in}n^{s+1}}\\
 & & \hspace{-0.2in} \mbox{} + O(n^s\log^{\delta(s+1,k)-1}\hspace{-0.05in}n)+O(\theta(u_1)n^s\log^{\delta(s+1,k)}\hspace{-0.05in}n)+O(\theta(u_1)n^s\log^{\delta(s,k)+1}\hspace{-0.05in}n)\\
 &= & A_{s+1,k}\prod_{i=1}^{k-1}f_{s+1,k,i}(u_i)n^{s+1}+O(\theta(u_1)n^s\log^{\delta(s+1,k)}\hspace{-0.04in}n)
\end{eqnarray*}
by a simple computation, which shows that the formula is true for $s+1$ and we complete the proof.\vspace{0.2in}\\
{\bf Acknowledgement.} I would like to thank Professor L\'{a}szl\'{o} T\'{o}th for helpful discussions of the topics of this paper.

{ \footnotesize {\sc University of Houston-Victoria, Department of Mathematics and Computer\\ Science,
 School of Arts \& Sciences, 14000 University Blvd, Sugar Land,TX, 77479 }

  {\em E-mail:}
{\bf huj@uhv.edu   }}

\end{document}